\documentclass[10pt]{article}
\usepackage{graphicx,amsthm,enumerate}
\usepackage{verbatim,amsmath,amscd,amssymb,color}
\usepackage{amsmath,amssymb,amsthm,color,eepic,eclarith}
\theoremstyle{plain}
\newcommand{\Pic}{\mathbf{P}(\mu, \nu \setminus \lambda)}
\newcommand{\PicA}{\mathbf{P}(\mu, \nu \setminus \lambda;A,A')}
\newcommand{\Cry}{\mathbf{B}(\mu)^{\nu}_{\lambda}}
\newcommand{\CryA}{\mathbf{B}(\mu)^{\nu}_{\lambda}[A']}
\newcommand{\MECry}{\fsquare(0.3cm,i_1) \otimes \fsquare(0.3cm,i_2) \otimes \cdots \otimes  \fsquare(0.3cm,i_k) \otimes \cdots \otimes \fsquare(0.3cm,i_N)}

\font\germ=eufm10

\usepackage{graphicx,amsthm}
\usepackage{verbatim,amsmath,amscd,amssymb,color}

\begin{document}
\font\germ=eufm10
\def\bl{\bullet}
\def\aaa{@}
\title{\Large\bf Admissible Pictures and 
Littlewood-Richardson Crystals}
\vskip0.6cm
\author{Toshiki N\textsc{akashima}
\thanks{Department of Mathematics, 
Sophia University, Kioicho 7-1, Chiyoda-ku,
Tokyo 102-8554, Japan.\hfill\break
\qquad E-mail: toshiki{\aaa}mm.sophia.ac.jp\quad
:supported in part by  JSPS Grants in 
Aid for Scientific Research \#19540050.
}
\and Miki S\textsc{himojo}
\thanks{Department of Mathematics, 
Sophia University, Kioicho 7-1, Chiyoda-ku,
Tokyo 102-8554, Japan.\hfill\break
E-mail: m-shimoj{\aaa}sophia.ac.jp}}

\date{}
\maketitle

\centerline{\it Dedicated to Professor Tetsuji Miwa 
on the occasion of his 60th birthday}

\abstract{We present a one-to-one correspondence between 
the set of admissible pictures and the Littlewood-Richardson 
crystals. As a simple consequence, we shall show that 
the set of pictures does not depend on the choice of 
admissible orders.}

    


\newcommand{\cI}{{\mathcal I}}
\newcommand{\cA}{{\mathcal A}}
\newcommand{\cB}{{\mathcal B}}
\newcommand{\cC}{{\mathcal C}}
\newcommand{\cD}{{\mathcal D}}
\newcommand{\cF}{{\mathcal F}}
\newcommand{\cH}{{\mathcal H}}
\newcommand{\cK}{{\mathcal K}}
\newcommand{\cL}{{\mathcal L}}
\newcommand{\cM}{{\mathcal M}}
\newcommand{\cN}{{\mathcal N}}
\newcommand{\cO}{{\mathcal O}}
\newcommand{\cS}{{\mathcal S}}
\newcommand{\cV}{{\mathcal V}}
\newcommand{\fra}{\mathfrak a}
\newcommand{\frb}{\mathfrak b}
\newcommand{\frc}{\mathfrak c}
\newcommand{\frd}{\mathfrak d}
\newcommand{\fre}{\mathfrak e}
\newcommand{\frf}{\mathfrak f}
\newcommand{\frg}{\mathfrak g}
\newcommand{\frh}{\mathfrak h}
\newcommand{\fri}{\mathfrak i}
\newcommand{\frj}{\mathfrak j}
\newcommand{\frk}{\mathfrak k}
\newcommand{\frI}{\mathfrak I}
\newcommand{\fm}{\mathfrak m}
\newcommand{\frn}{\mathfrak n}
\newcommand{\frp}{\mathfrak p}
\newcommand{\fq}{\mathfrak q}
\newcommand{\frr}{\mathfrak r}
\newcommand{\frs}{\mathfrak s}
\newcommand{\frt}{\mathfrak t}
\newcommand{\fru}{\mathfrak u}
\newcommand{\frA}{\mathfrak A}
\newcommand{\frB}{\mathfrak B}
\newcommand{\frF}{\mathfrak F}
\newcommand{\frG}{\mathfrak G}
\newcommand{\frH}{\mathfrak H}
\newcommand{\frJ}{\mathfrak J}
\newcommand{\frN}{\mathfrak N}
\newcommand{\frP}{\mathfrak P}
\newcommand{\frT}{\mathfrak T}
\newcommand{\frU}{\mathfrak U}
\newcommand{\frV}{\mathfrak V}
\newcommand{\frX}{\mathfrak X}
\newcommand{\frY}{\mathfrak Y}
\newcommand{\frZ}{\mathfrak Z}
\newcommand{\rA}{\mathrm{A}}
\newcommand{\rC}{\mathrm{C}}
\newcommand{\rd}{\mathrm{d}}
\newcommand{\rB}{\mathrm{B}}
\newcommand{\rD}{\mathrm{D}}
\newcommand{\rE}{\mathrm{E}}
\newcommand{\rH}{\mathrm{H}}
\newcommand{\rK}{\mathrm{K}}
\newcommand{\rL}{\mathrm{L}}
\newcommand{\rM}{\mathrm{M}}
\newcommand{\rN}{\mathrm{N}}
\newcommand{\rR}{\mathrm{R}}
\newcommand{\rT}{\mathrm{T}}
\newcommand{\rZ}{\mathrm{Z}}
\newcommand{\bbA}{\mathbb A}
\newcommand{\bbC}{\mathbb C}
\newcommand{\bbG}{\mathbb G}
\newcommand{\bbF}{\mathbb F}
\newcommand{\bbH}{\mathbb H}
\newcommand{\bbP}{\mathbb P}
\newcommand{\bbN}{\mathbb N}
\newcommand{\bbQ}{\mathbb Q}
\newcommand{\bbR}{\mathbb R}
\newcommand{\bbV}{\mathbb V}
\newcommand{\bbZ}{\mathbb Z}
\newcommand{\adj}{\operatorname{adj}}
\newcommand{\Ad}{\mathrm{Ad}}
\newcommand{\Ann}{\mathrm{Ann}}
\newcommand{\rcris}{\mathrm{cris}}
\newcommand{\ch}{\mathrm{ch}}
\newcommand{\coker}{\mathrm{coker}}
\newcommand{\diag}{\mathrm{diag}}
\newcommand{\Diff}{\mathrm{Diff}}
\newcommand{\Dist}{\mathrm{Dist}}
\newcommand{\rDR}{\mathrm{DR}}
\newcommand{\ev}{\mathrm{ev}}
\newcommand{\Ext}{\mathrm{Ext}}
\newcommand{\cExt}{\mathcal{E}xt}
\newcommand{\fin}{\mathrm{fin}}
\newcommand{\Frac}{\mathrm{Frac}}
\newcommand{\GL}{\mathrm{GL}}
\newcommand{\Hom}{\mathrm{Hom}}
\newcommand{\hd}{\mathrm{hd}}
\newcommand{\rht}{\mathrm{ht}}
\newcommand{\id}{\mathrm{id}}
\newcommand{\im}{\mathrm{im}}
\newcommand{\inc}{\mathrm{inc}}
\newcommand{\ind}{\mathrm{ind}}
\newcommand{\coind}{\mathrm{coind}}
\newcommand{\Lie}{\mathrm{Lie}}
\newcommand{\Max}{\mathrm{Max}}
\newcommand{\mult}{\mathrm{mult}}
\newcommand{\op}{\mathrm{op}}
\newcommand{\ord}{\mathrm{ord}}
\newcommand{\pt}{\mathrm{pt}}
\newcommand{\qt}{\mathrm{qt}}
\newcommand{\rad}{\mathrm{rad}}
\newcommand{\res}{\mathrm{res}}
\newcommand{\rgt}{\mathrm{rgt}}
\newcommand{\rk}{\mathrm{rk}}
\newcommand{\SL}{\mathrm{SL}}
\newcommand{\soc}{\mathrm{soc}}
\newcommand{\Spec}{\mathrm{Spec}}
\newcommand{\St}{\mathrm{St}}
\newcommand{\supp}{\mathrm{supp}}
\newcommand{\Tor}{\mathrm{Tor}}
\newcommand{\Tr}{\mathrm{Tr}}
\newcommand{\wt}{\mathrm{wt}}
\newcommand{\Ab}{\mathbf{Ab}}
\newcommand{\Alg}{\mathbf{Alg}}
\newcommand{\Grp}{\mathbf{Grp}}
\newcommand{\Mod}{\mathbf{Mod}}
\newcommand{\Sch}{\mathbf{Sch}}\newcommand{\bfmod}{{\bf mod}}
\newcommand{\Qc}{\mathbf{Qc}}
\newcommand{\Rng}{\mathbf{Rng}}
\newcommand{\Top}{\mathbf{Top}}
\newcommand{\Var}{\mathbf{Var}}
\newcommand{\BB}{\mathbf{B}}
\newcommand{\gromega}{\langle\omega\rangle}
\newcommand{\lbr}{\begin{bmatrix}}
\newcommand{\rbr}{\end{bmatrix}}
\newcommand{\forb}{\bigcirc\kern-2.8ex \because}
\newcommand{\forbb}{\bigcirc\kern-3.0ex \because}
\newcommand{\forbbb}{\bigcirc\kern-3.1ex \because}
\newcommand{\cd}{commutative diagram }
\newcommand{\SpS}{spectral sequence}
\newcommand\C{\mathbb C}
\newcommand\hh{{\hat{H}}}
\newcommand\eh{{\hat{E}}}
\newcommand\F{\mathbb F}
\newcommand\fh{{\hat{F}}}

\def\AA{{\mathcal A}}
\def\al{\alpha}
\def\bq{B_q(\ge)}
\def\bqm{B_q^-(\ge)}
\def\bqz{B_q^0(\ge)}
\def\bqp{B_q^+(\ge)}
\def\beneme{\begin{enumerate}}
\def\beq{\begin{equation}}
\def\beqn{\begin{eqnarray}}
\def\beqnn{\begin{eqnarray*}}
\def\bigsl{{\hbox{\fontD \char'54}}}
\def\bbra#1,#2,#3{\left\{\begin{array}{c}\hspace{-5pt}
#1;#2\\ \hspace{-5pt}#3\end{array}\hspace{-5pt}\right\}}
\def\cd{\cdots}
\def\CC{\hbox{\bf C}}
\def\ddd{\hbox{\germ D}}
\def\del{\delta}
\def\Del{\Delta}
\def\Delr{\Delta^{(r)}}
\def\Dell{\Delta^{(l)}}
\def\Delb{\Delta^{(b)}}
\def\Deli{\Delta^{(i)}}
\def\Delre{\Delta^{\rm re}}
\def\ei{e_i}
\def\eit{\tilde{e}_i}
\def\eneme{\end{enumerate}}
\def\ep{\epsilon}
\def\eeq{\end{equation}}
\def\eeqn{\end{eqnarray}}
\def\eeqnn{\end{eqnarray*}}
\def\fit{\tilde{f}_i}
\def\FF{{\rm F}}
\def\ft{\tilde{f}}
\def\gau#1,#2{\left[\begin{array}{c}\hspace{-5pt}#1\\
\hspace{-5pt}#2\end{array}\hspace{-5pt}\right]}
\def\ge{\hbox{\germ g}}
\def\gl{\hbox{\germ gl}}
\def\hom{{\hbox{Hom}}}
\def\ify{\infty}
\def\io{\iota}
\def\kp{k^{(+)}}
\def\km{k^{(-)}}
\def\llra{\relbar\joinrel\relbar\joinrel\relbar\joinrel\rightarrow}
\def\lan{\langle}
\def\lar{\longrightarrow}
\def\max{{\rm max}}
\def\ME{\rm ME}
\def\lm{\lambda}
\def\Lm{\Lambda}
\def\mapright#1{\smash{\mathop{\longrightarrow}\limits^{#1}}}
\def\mm{{\bf{\rm m}}}
\def\nd{\noindent}
\def\nn{\nonumber}
\def\nnn{\hbox{\germ n}}
\def\catob{{\mathcal O}(B)}
\def\oint{{\mathcal O}_{\rm int}(\ge)}
\def\ot{\otimes}
\def\op{\oplus}
\def\opi{\ovl\pi_{\lm}}
\def\ovl{\overline}
\def\plm{\Psi^{(\lm)}_{\io}}
\def\qq{\qquad}
\def\q{\quad}
\def\qed{\hfill\framebox[2mm]{}}
\def\QQ{\hbox{\bf Q}}
\def\qi{q_i}
\def\qii{q_i^{-1}}
\def\ra{\rightarrow}
\def\ran{\rangle}
\def\rlm{r_{\lm}}
\def\ssl{\mathfrak{sl}}
\def\slh{\widehat{\ssl_2}}
\def\ge{\hbox{\germ g}}
\def\ti{t_i}
\def\tii{t_i^{-1}}
\def\til{\tilde}
\def\tm{\times}
\def\tt{{\hbox{\germ{t}}}}
\def\ttt{\hbox{\germ t}}
\def\ua{U_{\AA}}
\def\ue{U_{\vep}}
\def\uq{U_q(\ge)}
\def\ufin{U^{\rm fin}_{\vep}}
\def\ufinp{(U^{\rm fin}_{\vep})^+}
\def\ufinm{(U^{\rm fin}_{\vep})^-}
\def\ufinz{(U^{\rm fin}_{\vep})^0}
\def\uqm{U^-_q(\ge)}
\def\uqp{U^+_q(\ge)}
\def\uqmq{{U^-_q(\ge)}_{\bf Q}}
\def\uqpm{U^{\pm}_q(\ge)}
\def\uqq{U_{\bf Q}^-(\ge)}
\def\uqz{U^-_{\bf Z}(\ge)}
\def\ures{U^{\rm res}_{\AA}}
\def\urese{U^{\rm res}_{\vep}}
\def\uresez{U^{\rm res}_{\vep,\ZZ}}
\def\util{\widetilde\uq}
\def\uup{U^{\geq}}
\def\ulow{U^{\leq}}
\def\bup{B^{\geq}}
\def\blow{\ovl B^{\leq}}
\def\vep{\varepsilon}
\def\vp{\varphi}
\def\vpi{\varphi^{-1}}
\def\VV{{\mathcal V}}
\def\xii{\xi^{(i)}}
\def\Xiioi{\Xi_{\io}^{(i)}}
\def\WW{{\mathcal W}}
\def\wtil{\widetilde}
\def\what{\widehat}
\def\wpi{\widehat\pi_{\lm}}
\def\ZZ{\mathbb Z}
\def\spsp(#1,#2){\begin{pmatrix}
#1,\\ \hline#2\end{pmatrix}}

\theoremstyle{definition}
\newtheorem{df}{\bf Definition}[section]
\newtheorem{pro}[df]{\bf Proposition}
\newtheorem{thm}[df]{\bf Theorem}
\newtheorem{lem}[df]{\bf Lemma}
\newtheorem{ex}[df]{\bf Example}
\newtheorem{cor}[df]{\bf Corollary}
\newtheorem{con}[df]{Conjecture}

\def\m@th{\mathsurround=0pt}

\def\fsquare(#1,#2){
\hbox{\vrule$\hskip-0.4pt\vcenter to #1{\normalbaselines\m@th
\hrule\vfil\hbox to #1{\hfill$\scriptstyle #2$\hfill}\vfil\hrule}$\hskip-0.4pt
\vrule}}

\def\addsquare(#1,#2){\hbox{$
	\dimen1=#1 \advance\dimen1 by -0.8pt
	\vcenter to #1{\hrule height0.4pt depth0.0pt%
	\hbox to #1{%
	\vbox to \dimen1{\vss%
	\hbox to \dimen1{\hss$\scriptstyle~#2~$\hss}%
	\vss}%
	\vrule width0.4pt}%
	\hrule height0.4pt depth0.0pt}$}}

\def\Fsquare(#1,#2){
\hbox{\vrule$\hskip-0.4pt\vcenter to #1{\normalbaselines\m@th
\hrule\vfil\hbox to #1{\hfill$#2$\hfill}\vfil\hrule}$\hskip-0.4pt
\vrule}}

\def\Addsquare(#1,#2){\hbox{$
	\dimen1=#1 \advance\dimen1 by -0.8pt
	\vcenter to #1{\hrule height0.4pt depth0.0pt%
	\hbox to #1{%
	\vbox to \dimen1{\vss%
	\hbox to \dimen1{\hss$~#2~$\hss}%
	\vss}%
	\vrule width0.4pt}%
	\hrule height0.4pt depth0.0pt}$}}

\def\hfourbox(#1,#2,#3,#4){%
	\fsquare(0.3cm,#1)\addsquare(0.3cm,#2)\addsquare(0.3cm,#3)\addsquare(0.3cm,#4)}

\def\Hfourbox(#1,#2,#3,#4){%
	\Fsquare(0.4cm,#1)\Addsquare(0.4cm,#2)\Addsquare(0.4cm,#3)\Addsquare(0.4cm,#4)}

\def\HHfourbox(#1,#2,#3,#4){%
	\Fsquare(0.8cm,#1)\Addsquare(0.8cm,#2)\Addsquare(0.8cm,#3)\Addsquare(0.8cm,#4)}

\def\fsq(#1){%
          \fsquare(0.3cm,#1)}

\def\hthreebox(#1,#2,#3){%
	\fsquare(0.3cm,#1)\addsquare(0.3cm,#2)\addsquare(0.3cm,#3)}

\def\htwobox(#1,#2){%
	\fsquare(0.3cm,#1)\addsquare(0.3cm,#2)}

\def\vfourbox(#1,#2,#3,#4){%
	\normalbaselines\m@th\offinterlineskip
	\vcenter{\hbox{\fsquare(0.3cm,#1)}
	      \vskip-0.4pt
	      \hbox{\fsquare(0.3cm,#2)}	
	      \vskip-0.4pt
	      \hbox{\fsquare(0.3cm,#3)}	
	      \vskip-0.4pt
	      \hbox{\fsquare(0.3cm,#4)}}}

\def\VVfourbox(#1,#2,#3,#4){%
	\normalbaselines\m@th\offinterlineskip
	\vcenter{\hbox{\Fsquare(0.8cm,#1)}
	      \vskip-0.4pt
	      \hbox{\Fsquare(0.8cm,#2)}	
	      \vskip-0.4pt
	      \hbox{\Fsquare(0.8cm,#3)}	
	      \vskip-0.4pt
	      \hbox{\Fsquare(0.8cm,#4)}}}

\def\Vfourbox(#1,#2,#3,#4){%
	\normalbaselines\m@th\offinterlineskip
	\vcenter{\hbox{\Fsquare(0.4cm,#1)}
	      \vskip-0.4pt
	      \hbox{\Fsquare(0.4cm,#2)}	
	      \vskip-0.4pt
	      \hbox{\Fsquare(0.4cm,#3)}	
	      \vskip-0.4pt
	      \hbox{\Fsquare(0.4cm,#4)}}}

\def\vthreebox(#1,#2,#3){%
	\normalbaselines\m@th\offinterlineskip
	\vcenter{\hbox{\fsquare(0.3cm,#1)}
	      \vskip-0.4pt
	      \hbox{\fsquare(0.3cm,#2)}	
	      \vskip-0.4pt
	      \hbox{\fsquare(0.3cm,#3)}}}

\def\vtwobox(#1,#2){%
	\normalbaselines\m@th\offinterlineskip
	\vcenter{\sbox{\fsquare(0.3cm,#1)}
	      \vskip-0.4pt
	      \hbox{\fsquare(0.3cm,#2)}}}

\def\vtwobox2(#1,#2){%
	\normalbaselines\m@th\offinterlineskip
	\vcenter{\hbox{#1}
	      \vskip-0.4pt
	      \hbox{\fsquare(0.3cm,#2)}}}

\def\Hthreebox(#1,#2,#3){%
	\Fsquare(0.4cm,#1)\Addsquare(0.4cm,#2)\Addsquare(0.4cm,#3)}

\def\HHthreebox(#1,#2,#3){%
	\Fsquare(0.8cm,#1)\Addsquare(0.8cm,#2)\Addsquare(0.8cm,#3)}

\def\Htwobox(#1,#2){%
	\Fsquare(0.4cm,#1)\Addsquare(0.4cm,#2)}

\def\H6twobox(#1,#2){%
	\Fsquare(0.6cm,#1)\Addsquare(0.6cm,#2)}

\def\HHtwobox(#1,#2){%
	\Fsquare(0.8cm,#1)\Addsquare(0.8cm,#2)}

\def\Vthreebox(#1,#2,#3){%
	\normalbaselines\m@th\offinterlineskip
	\vcenter{\hbox{\Fsquare(0.4cm,#1)}
	      \vskip-0.4pt
	      \hbox{\Fsquare(0.4cm,#2)}	
	      \vskip-0.4pt
	      \hbox{\Fsquare(0.4cm,#3)}}}

\def\Vtwobox(#1,#2){%
	\normalbaselines\m@th\offinterlineskip
	\vcenter{\hbox{\Fsquare(0.4cm,#1)}
	      \vskip-0.4pt
	      \hbox{\Fsquare(0.4cm,#2)}}}

\def\twoone(#1,#2,#3){%
	\normalbaselines\m@th\offinterlineskip
	\vcenter{\hbox{\htwobox({#1},{#2})}
	      \vskip-0.4pt
	      \hbox{\fsquare(0.3cm,#3)}}}

\def\twothree(#1,#2,#3,#4,#5){%
	\normalbaselines\m@th\offinterlineskip
	\vcenter{\hbox{\htwobox({#1},{#2})}
	      \vskip-0.4pt
	      \hbox{\hthreebox({#3},{#4},{#5})}}}

\def\threethree(#1,#2,#3,#4,#5,#6){%
	\normalbaselines\m@th\offinterlineskip
	\vcenter{\hbox{\hthreebox({#1},{#2},{#3})}
	      \vskip-0.4pt
	      \hbox{\hthreebox({#4},{#5},{#6})}}}

\def\threethreeone(#1,#2,#3,#4,#5,#6,#7){%
	\normalbaselines\m@th\offinterlineskip
	\vcenter{\hbox{\hthreebox({#1},{#2},{#3})}
	      \vskip-0.4pt
	      \hbox{\hthreebox({#4},{#5},{#6})}
              \vskip-0.4pt
              \hbox{\fsquare(0.3cm,#7)}}}

\def\threethreetwo(#1,#2,#3,#4,#5,#6,#7,#8){%
	\normalbaselines\m@th\offinterlineskip
	\vcenter{\hbox{\hthreebox({#1},{#2},{#3})}
	      \vskip-0.4pt
	      \hbox{\hthreebox({#4},{#5},{#6})}
              \vskip-0.4pt
              \hbox{\htwobox({#7},{#8})}}}

\def\Twoone(#1,#2,#3){%
	\normalbaselines\m@th\offinterlineskip
	\vcenter{\hbox{\Htwobox({#1},{#2})}
	      \vskip-0.4pt
	      \hbox{\Fsquare(0.4cm,#3)}}}

\def\TTwoone(#1,#2,#3){%
	\normalbaselines\m@th\offinterlineskip
	\vcenter{\hbox{\H6twobox({#1},{#2})}
	      \vskip-0.4pt
	      \hbox{\Fsquare(0.6cm,#3)}}}

\def\threeone(#1,#2,#3,#4){%
	\normalbaselines\m@th\offinterlineskip
	\vcenter{\hbox{\hthreebox({#1},{#2},{#3})}
	      \vskip-0.4pt
	      \hbox{\fsquare(0.3cm,#4)}}}

\def\Threeone(#1,#2,#3,#4){%
	\normalbaselines\m@th\offinterlineskip
	\vcenter{\hbox{\Hthreebox({#1},{#2},{#3})}
	      \vskip-0.4pt
	      \hbox{\Fsquare(0.4cm,#4)}}}

\def\Threetwo(#1,#2,#3,#4,#5){%
	\normalbaselines\m@th\offinterlineskip
	\vcenter{\hbox{\Hthreebox({#1},{#2},{#3})}
	      \vskip-0.4pt
	      \hbox{\Htwobox({#4},{#5})}}}

\def\threetwo(#1,#2,#3,#4,#5){%
	\normalbaselines\m@th\offinterlineskip
	\vcenter{\hbox{\hthreebox({#1},{#2},{#3})}
	      \vskip-0.4pt
	      \hbox{\htwobox({#4},{#5})}}}

\def\twotwo(#1,#2,#3,#4){%
	\normalbaselines\m@th\offinterlineskip
	\vcenter{\hbox{\htwobox({#1},{#2})}
	      \vskip-0.4pt
	      \hbox{\htwobox({#3},{#4})}}}

\def\Twotwo(#1,#2,#3,#4){%
	\normalbaselines\m@th\offinterlineskip
	\vcenter{\hbox{\Htwobox({#1},{#2})}
	      \vskip-0.4pt
	      \hbox{\Htwobox({#3},{#4})}}}

\def\TTwotwo(#1,#2,#3,#4){%
	\normalbaselines\m@th\offinterlineskip
	\vcenter{\hbox{\H6twobox({#1},{#2})}
	      \vskip-0.4pt
	      \hbox{\H6twobox({#3},{#4})}}}

\def\twooneone(#1,#2,#3,#4){%
	\normalbaselines\m@th\offinterlineskip
	\vcenter{\hbox{\htwobox({#1},{#2})}
	      \vskip-0.4pt
	      \hbox{\fsquare(0.3cm,#3)}
	      \vskip-0.4pt
	      \hbox{\fsquare(0.3cm,#4)}}}

\def\Twooneone(#1,#2,#3,#4){%
	\normalbaselines\m@th\offinterlineskip
	\vcenter{\hbox{\Htwobox({#1},{#2})}
	      \vskip-0.4pt
	      \hbox{\Fsquare(0.4cm,#3)}
	      \vskip-0.4pt
	      \hbox{\Fsquare(0.4cm,#4)}}}

\def\Twotwoone(#1,#2,#3,#4,#5){%
	\normalbaselines\m@th\offinterlineskip
	\vcenter{\hbox{\Htwobox({#1},{#2})}
	      \vskip-0.4pt
	      \hbox{\Htwobox({#3},{#4})}
              \vskip-0.4pt
	      \hbox{\Fsquare(0.4cm,#5)}}}

\def\twotwoone(#1,#2,#3,#4,#5){%
	\normalbaselines\m@th\offinterlineskip
	\vcenter{\hbox{\htwobox({#1},{#2})}
	      \vskip-0.4pt
	      \hbox{\htwobox({#3},{#4})}
              \vskip-0.4pt
	      \hbox{\fsquare(0.3cm,#5)}}}

\def\twotwotwo(#1,#2,#3,#4,#5,#6){%
	\normalbaselines\m@th\offinterlineskip
	\vcenter{\hbox{\htwobox({#1},{#2})}
	      \vskip-0.4pt
	      \hbox{\htwobox({#3},{#4})}
              \vskip-0.4pt
	      \hbox{\htwobox({#5},{#6})}}}

\def\threetwoone(#1,#2,#3,#4,#5,#6){%
	\normalbaselines\m@th\offinterlineskip
	\vcenter{\hbox{\hthreebox({#1},{#2},{#3})}
	      \vskip-0.4pt
	      \hbox{\htwobox({#4},{#5})}
              \vskip-0.4pt
	      \hbox{\fsquare(0.3cm,#6)}}}

\def\fourthreetwoone{%
	\normalbaselines\m@th\offinterlineskip
	\vcenter{\hbox{\hfourbox(,,,)}
	      \vskip-0.4pt
              \hbox{\hthreebox(,,)}
              \vskip-0.4pt
	      \hbox{\htwobox(,)}
              \vskip-0.4pt
	      \hbox{\fsquare(0.3cm,)}}}

\def\fourtwoone(#1,#2,#3,#4,#5,#6,#7){%
	\normalbaselines\m@th\offinterlineskip
	\vcenter{\hbox{\hfourbox({#1},{#2},{#3},{#4})}
	      \vskip-0.4pt
	      \hbox{\htwobox({#5},{#6})}
              \vskip-0.4pt
	      \hbox{\fsquare(0.3cm,{#7})}}}

\def\fourthreetwo(#1,#2,#3,#4,#5,#6,#7,#8,#9){%
	\normalbaselines\m@th\offinterlineskip
	\vcenter{\hbox{\hfourbox({#1},{#2},{#3},{#4})}
	      \vskip-0.4pt
              \hbox{\hthreebox({#5},{#6},{#7})}
              \vskip-0.4pt
	      \hbox{\htwobox({#8},{#9})}}}
\def\fourthreeone(#1,#2,#3,#4,#5,#6,#7,#8){%
	\normalbaselines\m@th\offinterlineskip
	\vcenter{\hbox{\hfourbox({#1},{#2},{#3},{#4})}
	      \vskip-0.4pt
              \hbox{\hthreebox({#5},{#6},{#7})}
              \vskip-0.4pt
	      \hbox{\fsquare(0.3cm,{#8})}}}

\def\TThreetwoone(#1,#2,#3,#4,#5,#6){%
	\normalbaselines\m@th\offinterlineskip
	\vcenter{\hbox{\HHthreebox({#1},{#2},{#3})}
	      \vskip-0.4pt
	      \hbox{\HHtwobox({#4},{#5})}
              \vskip-0.4pt
	      \hbox{\fsquare(0.8cm,#6)}}}

\def\Threeoneone(#1,#2,#3,#4,#5){%
	\normalbaselines\m@th\offinterlineskip
	\vcenter{\hbox{\Hthreebox({#1},{#2},{#3})}
	      \vskip-0.4pt
	      \hbox{\Fsquare(0.4cm,#4)}
              \vskip-0.4pt
	      \hbox{\Fsquare(0.4cm,#5)}}}

\def\threeoneone(#1,#2,#3,#4,#5){%
	\normalbaselines\m@th\offinterlineskip
	\vcenter{\hbox{\hthreebox({#1},{#2},{#3})}
	      \vskip-0.4pt
	      \hbox{\fsquare(0.3cm,#4)}
              \vskip-0.4pt
	      \hbox{\fsquare(0.3cm,#5)}}}

\def\FFourtwoone(#1,#2,#3,#4,#5,#6,#7){%
	\normalbaselines\m@th\offinterlineskip
	\vcenter{\hbox{\HHfourbox({#1},{#2},{#3},{#4})}
	      \vskip-0.4pt
	      \hbox{\HHtwobox({#5},{#6})}
              \vskip-0.4pt
	      \hbox{\fsquare(0.8cm,#7)}}}

\def\a{\fsquare(0.3cm){1}\addsquare(0.3cm)(2)\addsquare(0.3cm)(3)}

\def\b{\hbox{%
	\normalbaselines\m@th\offinterlineskip
	\vcenter{\hbox{\fsquare(0.3cm){2}}\vskip-0.4pt\hbox{\fsquare(0.3cm){2}}}}}

\def\c{\hbox{\normalbaselines\m@th\offinterlineskip%
	\vcenter{\hbox{\a}\vskip-0.4pt\hbox{\b}}}}


\dimen1=0.4cm\advance\dimen1 by -0.8pt
\def\ffsquare#1{%
	\fsquare(0.4cm,\hbox{#1})}

\def\naga{%
	\hbox{$\vcenter to 0.4cm{\normalbaselines\m@th
	\hrule\vfil\hbox to 1.2cm{\hfill$\cdots$\hfill}\vfil\hrule}$}}

\def\vnaga{\normalbaselines\m@th\baselineskip0pt\offinterlineskip%
	\vrule\vbox to 1.2cm{\vskip7pt\hbox to \dimen1{$\hfil\vdots\hfil$}\vfil}\vrule}

\def\dhbox{\hbox{$\ffsquare 1 \naga \ffsquare N$}}

\def\dvbox{\hbox{\normalbaselines\m@th\baselineskip0pt\offinterlineskip\vbox{%
	  \hbox{$\ffsquare 1$}\vskip-0.4pt\hbox{$\vnaga$}\vskip-0.4pt\hbox{$\ffsquare N$}}}}

\def\sq(#1){\fsquare(0.4cm,#1)}
\def\Sq(#1){\fsquare(0.5cm,#1)}
\def\SSq(#1){\fsquare(0.9cm,#1)}
\def\Sqj{\Sq(j)}
\def\Sqjb{\Sq(\ovl j)}
\def\Sqjm{\Sq(\scriptstyle{j-1})}
\def\Sqjp{\Sq(\scriptstyle{j+1})}
\def\Sqjmb{\Sq(\scriptstyle{\ovl{j-1}})}
\def\SSqj{\SSq(j)}
\def\Sqjpb{\Sq(\scriptstyle{\ovl{j+1}})}
\def\SSqjb{\SSq(\ovl j)}
\def\SSqjm{\SSq(\scriptstyle{j-1})}
\def\SSqjp{\SSq(\scriptstyle{j+1})}
\def\SSqjmb{\SSq(\scriptstyle{\ovl{j-1}})}
\def\SSqjpb{\SSq(\scriptstyle{\ovl{j+1}})}

\def\mapright#1{\smash{\mathop{\longrightarrow}\limits^{#1}}}
\def\map#1{\smash{\mathop{\longmapsto}\limits^{#1}}}

\renewcommand{\thesection}{\arabic{section}}
\section{Introduction}
\setcounter{equation}{0}
\renewcommand{\theequation}{\thesection.\arabic{equation}}

On (skew)Young diagrams, we introduce several orders,
{\it e.g.,} $\leqslant_P$, $\leqslant_J$, $\leqslant_A$, 
etc (see \ref{pic}) to treat our main subject ''picture'', which 
is a bijective map between two skew Young diagrams,
which preserves the order in the following sense:
$f:\leqslant_P\to\leqslant_J$ and $f^{-1}:\leqslant_P\to\leqslant_J$ (\cite{CS},\cite{JP},\cite{Z}).
For Young diagrams $\lm,\mu,\nu$ with $|\lm|+|\mu|=|\nu|$, 
let $\Pic$ be the set of pictures from $\mu$ 
to $\nu\setminus\lm$ and $\Cry$ the 
Littlewood-Richardson crystal as in \cite{NS}
(see also \ref{subsec:LR}). 
Let $c_{\lm,\mu}^\nu$ be the usual 
Littlewood-Richardson number. Then, by the fact 
\[
 \sharp\Pic=c_{\lm,\mu}^\nu=\sharp\Cry,
\]
we deduced that there exists a bijection between 
$\Pic$ and $\Cry$.
It has been 
revealed in \cite{NS} that 
there exists a natural one-to-one correspondence between
$\Pic$ and $\Cry$.

We try to generalize the notion of pictures by using 
''admissible order'', which is an order in 
a certain class of 
total orders on a skew diagram
(or more generally, a subset of $\bbN\times\bbN$).
Indeed, the order $\leqslant_J$ is 
a sort of admissible orders.
In the last section of \cite{NS}, we define the new 
set of ``admissible pictures'' associated with admissible
orders $A$ on $\nu\setminus\lm$ and $A'$ on 
$\mu$, denoted by $\PicA$. 
We also get the Littilwood-Richardson crystal associated with 
an admissible order, denoted by $\CryA$. Then we conjectured:
\begin{con}[\cite{NS}]
Let $A$ (resp. $A'$) be an admissible order on 
$\nu\setminus\lm$ (resp. $\mu$). 
There exists a bijection 
\[
 \Psi:\mathbf{B}(\mu)_\lm^\nu[A']\longrightarrow
{\bf P}(\mu,\nu\setminus\lm:A,A').
\]
\end{con} 
The affirmative answer for this conjecture is 
given as Theorem \ref{main} in Sect.4 below.

In \cite{HK}, it has been shown 
that the Littlewood-Richardson 
crystal does not depend on the choice of admissible orders.
Furthermore, so does not the definition of the bijection $\Psi$.
Therefore, we obtain:
\begin{cor}
For arbitrary admissible orders $A$ on $\nu\setminus\lm$ and 
$A'$ on $\mu$,
\[
\Pic= {\bf P}(\mu,\nu\setminus\lm:A,A').
\]
\end{cor}
This result has already been obtained 
in \cite{CS2} and \cite{FG} for 
more general setting.
They used some purely combinatorial methods different from 
ours. Here it can be said that we give 
a new proof of the pictures' independence of admissible 
orders. Our main tool is a procedure ''addition'' 
obtained from the tensor products of crystals.  It plays a
crucial role in the proof, which realizes
the Littlewood-Richardson rules in terms of crystals and 
connects pictures and the Littlewood-Richardson crystals
directly.

We have obtained the Littlwood-Richardson crystals
for other classical types (\cite{N1}) 
in the similar description to the type $A_n$.
Hence, it allows us to expect that 
it is possible to generalize the notion ''pictures''
to other classical types.

The organization of the article is as follows.
In Sect.2, we prepare the ingredients treated in the paper, 
(skew) Young diagrams, Young tableaux, 
admissible orders and pictures. In Sect.3, we review 
the crystal-theoretical interpretation of Littlewood-Richardson
rules. The definitions of additions and readings are given.
The main theorem is stated in Sect.4.
The last three sections, 5,6 and 7 are 
devoted to show the main theorem.
\renewcommand{\thesection}{\arabic{section}}
\section{Pictures}
\setcounter{equation}{0}
\renewcommand{\theequation}{\thesection.\arabic{equation}}

\subsection{Young diagrams and Young tableaux}\label{pT}
Let $\lambda=(\lambda_1,\lambda_2,\cdots,\lambda_m)$
be a Young diagram or a partition, which satisfies
$\lambda_1 \geq\lambda_2 \geq \cdots \geq \lambda_m \geq 0$.
Let $\lambda$ and $\mu$ be Young diagrams with $\mu\subset\lm$.
A {\it skew diagram} $\lambda \setminus \mu$ is obtained by 
subtracting set-theoretically $\mu$ from $\lm$.

In this article we frequently consider 
a (skew) Young diagram as a subset of $\bbN\times\bbN$ by 
identifying the box in the $i$-th row and the $j$-th column 
with $(i,j)\in \bbN\times\bbN$.

\begin{ex}@A Young diagram $\lambda=(2,2,1)$ is expressed by 
$\left\{
(1,1), (1,2),
(2,1), (2,2),
(3,1)
\right\}$.
\end{ex}

As in \cite{F}, 
in the sequel,
a ''Young tableau'' means  a semi-standard tableau.
For a Young tableau $T$ of shape $\lm$, 
we also consider a ''coordinate '' in $\bbN\times\bbN$
like as $\lm$. 
Then an entry of $T$ in $(i,j)$ is denoted by $T_{i,j}$ and 
called $(i,j)$-entry.
For $k>0$, define (\cite{NS})
\begin{equation}
\label{Y-diag}
T^{(k)} = \{ (l,m) \in \lm | T_{l,m} = k \}.
\end{equation}
There is no two elements in one column in $T^{(k)}$. 
For a Young tableau $T$ with $(i,j)$-entry $T_{i,j}=k$, 
we define a function $p(T;i,j)$  (\cite{NS}) as
the number of $(i,j)$-entry 
from the right in $T^{(k)}$.
It is immediate from the definition:
\begin{equation}
\text{If }T_{i,j}=T_{x,y}\text{ and }p(T;i,j)=p(T;x,y),
\text{ then }(i,j)=(x,y)
\label{ijxy}
\end{equation}

\subsection{Picture}\label{pic}

First, we shall introduce the original notion of 
''picture'' as in \cite{Z}.

We define the following two kinds of orders on
a subset $X\subset\bbN \times \bbN$: For
$(a,b), \,\,(c,d)\in X$,
\begin{enumerate}
\item $(a,b)\leqslant_P(c,d)$ iff
$a\leq c \text{ and } b \leq d$.
\item $(a,b) \leqslant_J (c,d)$ iff
$a < c \text{, or }  a = c \text{ and } b\geq d$.
\end{enumerate}
Note that the order $\leqslant_P$ is a partial order 
and $\leqslant_J$ is a total order.
\begin{df}[\cite{Z}] 
Let $X,Y \subset \bbN \times \bbN$.
\begin{enumerate}
\item
A map $f : X \to Y$ is said to be 
{\it PJ-standard} if it satisfies
\[
\text{For }(a,b),(c,d) \in X, 
\text{if }(a,b) \leqslant_P (c,d), \text{ then }  
f(a,b) \leqslant_J f(c,d).
\]
\item
A map $f:X\to Y$ is a {\it picture} if it is 
bijective and
both $f$ and $f^{-1}$ are PJ-standard.
\end{enumerate}
\end{df}

Taking three Young diagrams
$\lambda, \mu, \nu\subset \bbN\times\bbN$, 
denote the set of pictures by:
\[
\Pic := \{ f : \mu \to \nu \setminus \lambda 
\,|\,f\text{ is a picture.}\}
\]

Next, we shall generalize the notion of pictures by using
a total order on a subset 
$X\subset \bbN\times \bbN$, called an
``{\it admissible order}'',
\begin{df}
\begin{enumerate}
\item
A total order $\leqslant_A$ on $X\subset\bbN\times\bbN$ 
is called {\it admissible}
if it satisfies:
\[
\text{For any }(a,b),\,\,(c,d)\in X\text{ if }
a\leq c\text{ and }b\geq d\text{ then } (a,b)\leqslant_A (c,d).
\]
\item
For $X,Y\subset\bbN\times\bbN$ and a map $f:X\to Y$, if $f$ satisfies
that if $(a,b)\leqslant_P(c,d)$, then $f(a,b)\leqslant_A f(c,d)$ for 
any $(a,b),~(c,d)\in X$, then $f$ is called $PA$-standard.
\item
Let $\leqslant_A$ (resp. $\leqslant_{A'}$) be an admissible order on 
$X\text{(resp. }Y)\subset\bbN\times\bbN$. A bijective map 
$f:X\to Y$ is called an $(A,A')$-{\it admissible picture} or 
simply, an {\it admissible picture} if $f$ is $PA$-standard and 
$f^{-1}$ is $PA'$-standard.
\end{enumerate}
\end{df}
{\sl Remark.} Note that for fixed $X\subset\bbN\times\bbN$, there 
can be several admissible orders on $X$. For example, 
the order $\leqslant_J$ is one of admissible orders on $X$. 
If we define the total order $\leqslant_F$ by 
\[
 (a,b)\leqslant_F (c,d) \text{ iff }b>d, 
\text{ or }b=d \text{ and }a<c, 
\]
which is also admissible.

For $X,\,Y\in\bbN\times\bbN$, let 
$\leqslant_A$ (resp. $\leqslant_{A'}$) be 
an admissible order on $X$
(resp. $Y$). 
We denote a set of $(A,A')$-admissible pictures by 
${\bf P}(X,Y:A,A')$.

\renewcommand{\thesection}{\arabic{section}}
\section{Crystals}
\setcounter{equation}{0}
\renewcommand{\theequation}{\thesection.\arabic{equation}}

The basic references for the theory of crystals are
\cite{K1},\cite{K2}.
\subsection{Readings and Additions}
Let $\BB=\{\fsq(i)\,|\,1\leq i\leq n+1\}$ 
be the crystal of the vector representation 
$V(\Lm_1)$ of the quantum group $U_q(A_n)$ (\cite{KN}).
As in \cite{NS}, we shall identify a dominant integral weight 
of type $A_n$ with a Young diagram in the standard way, 
{\it e.g.,} the fundamental weight $\Lm_1$ is identified with 
a square box $\fsq()$. For a Young diagram $\lm$, 
let $B(\lm)$ be the crystal of the finite-dimensional 
irreducible $U_q(A_n)$-module $V(\lm)$.
Set $N:=|\lm|$. Then there exists an embedding of crystals:
$B(\lm)\hookrightarrow \BB^N$ and an element in $B(\lm)$ 
is realized by a Young tableau of shape $\lm$
(\cite{KN}). Such an embedding is not unique. Indeed, 
they are called a 'reading' and described by:
\begin{df}[\cite{HK}]
Let $A$ be an admissible order on 
a Young diagram $\mu$ with $|\lm|=N$.
For $T\in B(\lm)$, by reading the entries 
in $T$ according to $A$, 
we obtain the map 
\[
 R_A:B(\lm)\longrightarrow B^{\ot N}\q (T\mapsto
\fsq(i_1)\ot\cd\ot\fsq(i_N))),
\]
which is called an {\it admissible reading}
associated with the order $A$.
The map $R_A$ is 
an embedding of crystals.
\end{df}
The following are typical readings.
\begin{df}
Let $T$ be an element in $B(\lm)$ of type $A_n$, namely, 
a Young tableau of shape $\lm$ 
with entries $\{1,2,\cd, n+1\}$.
\begin{enumerate}
\item
We read the entries in $T$ each row 
from right to left and from the top
row to the bottom row, that is, 
we read the entries according to the order $\leqslant_J$.
Then the resulting sequence of the entries 
$i_1,i_2,\cd,i_N$ gives the embedding of crystals:
\[
{\rm ME}(=R_J): B(\lm)\hookrightarrow \BB^{\ot N} \q
(T\mapsto \fsq(i_1)\ot\cd\ot\fsq(i_N)),
\]
which is called a {\it middle-eastern reading}.
\item
We read the entries in $T$ each column from the top to the 
bottom and from the right-most column to the left-most column, 
that is, 
we read the entries according to the order $\leqslant_F$.
Then  the resulting sequence of the entries 
$i_1,i_2,\cd,i_N$ gives the embedding of crystals:
\[
{\rm FE}(=R_F): B(\lm)\hookrightarrow \BB^{\ot N} \q
(T\mapsto \fsq(i_1)\ot\cd\ot\fsq(i_N)),
\]
which is called a {\it far-eastern reading}.
\end{enumerate}
\end{df}
\begin{df}
For $i\in\{1,2,\cd,n+1\}$ and a Young diagram 
$\lambda=(\lambda_1,\lambda_2,\cdots, \lambda_n)$, we define
\[
\lambda[i]:=
(\lambda_1,\lambda_2,\cdots,\lambda_i+1,\cdots, \lambda_n)
\]
which is said to be an {\it addition} of $i$ to $\lm$.
In general, 
for $i_1,i_2,\cd,i_N\in\{1,2,\cd,n+1\} $ and a Young diagram 
$\lm$, we define 
\[
\lambda[i_1,i_2,\cdots,i_N]
:=(\cd ((\lambda[i_1])[i_2])\cd)[i_N],
\]
which is called an {\it addition} of $i_1,\cd,i_N$ to $\lm$.
\end{df}

\begin{ex}
For a sequence ${\bf i}=31212$, 
the addition of $\bf i$ to $\lambda = \twoone(,,)$ is:
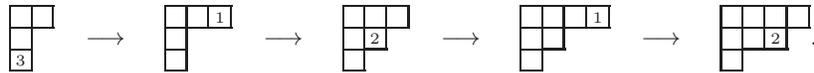
\[
\twooneone(,,,3) \quad \longrightarrow \quad \threeoneone(,,1,,)
\quad \longrightarrow \quad \threetwoone(,,,,2,)
\quad \longrightarrow \quad \fourtwoone(,,,1,,,) \quad \longrightarrow
 \quad \fourthreeone(,,,,,,2,).
\]
\end{ex}
\vspace{0.5cm}
{\sl Remark.}
For a Young diagram $\lm$, 
an addition $\lm[i_1,\cd,i_N]$ is not 
necessarily a Young diagram.
For instance, a sequence ${\bf i'}=22133$ and 
$\lambda=(2,2)$,
the addition $\lm[{\bf i'}]=(3,3,2)$ is a Young diagram.
But, in the second step of the addition, it becomes the diagram 
$\lm[2,2]=(2,3)$, which is not a Young diagram.

\subsection{Littlewood-Richardson Crystal}\label{subsec:LR}
As an application of the description of 
crystal bases of type $A_n$,
we see so-called ``Littlwood-Richardson rule'' of type $A_n$.

For a sequence ${\mathbf i}=
i_1,i_2,\cd,i_N\,\,(i_j\in \{1,2,\cd,n+1\})$ and 
a Young diagram $\lm$, let $\til\lm:=\lm[i_1,i_2,\cd,i_{N}]$ be
an addition of $i_1,i_2,\cd,i_N$ to $\lm$. Then set
\[
{\mathbf B}(\lm:{\mathbf i})=\begin{cases}
{\mathbf B}(\til\lm)&\text{ if }\lm[i_1,\cd,i_k]
\text{ is a Young diagram for any }k=1,2,\cd,N,\\
\emptyset&\text{otherwise.}
\end{cases}
\]

\begin{thm}[\cite{N1}]
\label{LRrule}
Let $\lambda$ and $\mu$ be Young diagrams with at most $n$ rows.
Then we have
\begin{equation}
\mathbf{B}(\lambda) \otimes \mathbf{B}(\mu) \cong  
\bigoplus_{\tiny\begin{array}{l}T\in \mathbf{B}(\mu),\\
{\rm FE}(T)= \fsq(i_1) 
\otimes\cdots\otimes\fsq(i_N)
\end{array}} 
\mathbf{B}(\lambda:i_1,i_2,\ldots,i_N).
\label{LR}
\end{equation}
\end{thm}
Here note that Theorem \ref{LRrule} is valid 
for an arbitrary admissible order
$A$, that is, in (\ref{LR}) we can replace ${\rm FE}(T)$ with 
$R_A(T)$. For an admissible order $A$ on $\mu$ define 
\[
\mathbf{B}(\mu)_\lm^\nu[A]
:= 
\left\{T\in{\mathbf B}(\mu)|
\begin{array}{l} 
R_A(T)=\MECry. \\
\text{For any }k=1,\cd,N,\\
\lm[i_1,\cd,i_k]\text{ is a Young diagram and}\\
\lm[i_1,\cd,i_N]=\nu.
\end{array}
\right\},
\]
Let $\Cry:=\Cry[J]$ as in \cite{NS}.
It is shown in \cite{HK} that for any admissible order 
$A$ on $\mu$,
\begin{equation}
\mathbf{B}(\mu)_\lm^\nu[A]=\Cry,
\label{crya=cry}
\end{equation}
which is called a {\it Littlewood-Richardson crystal}
associated with a triplet $(\lm,\mu,\nu)$.

\renewcommand{\thesection}{\arabic{section}}
\section{Main Theorem}
\setcounter{equation}{0}
\renewcommand{\theequation}{\thesection.\arabic{equation}}

For Young diagrams $\lambda, \mu, \nu$ with 
$|\lambda|+|\mu|=|\nu|$, we define the 
map $\Phi:\PicA\to\CryA$:
For $f=(f_1, f_2)\in\PicA$ where $f_1$ (resp. $f_2$) stands for 
the first (resp. second) coordinate in $\bbN\times\bbN$, set 
\[
\Phi(f)_{i,j}:= f_{1}(i,j),
\]
that is, $\Phi(f)$ is a filling of shape $\mu$ and 
its $(i, j)$-entry is given as $f_1(i, j)$.

Furthermore, for $T \in \CryA$, 
define a  map 
$\Psi : \CryA \to \PicA$ by 
\[
\Psi(T) : (i,j)\in\mu \mapsto (T_{i,j}, 
\lambda_{T_{i,j}}+p(T;i,j))\in\nu\setminus\lambda,
\]
where 
$p(T;i,j)$ as in \ref{pT}.

The following is the main result in this article, which 
is conjectured in \cite{NS}:
\begin{thm}\label{main}
For Young diagrams $\lm,\mu,\nu$ as above,
the maps
\[
 \Phi:\PicA\to\CryA,\qq
 \Psi:\CryA\to\PicA,
\] 
are bijections and 
they are inverse each other.
\end{thm}
It follows from  (\ref{crya=cry}) that the set 
$\CryA$ does not depend on the choice of an admissible order 
$A'$. Furthermore, it is easy to see from the definition 
that the map $\Psi$ does not depend on the choice of $A'$.
Therefore, we have 
\[
\Psi(\CryA)=\Psi(\Cry).
\]
By Theorem \ref{main} we have $\Psi(\CryA)=\PicA$ 
(resp. $\Psi(\Cry)=\Pic$). Hence, we obtain:
\begin{cor}
For arbitrary admissible orders $A$ on $\nu\setminus\lm$ and 
$A'$ on $\mu$,
\[
\PicA=\Pic.
\]
\end{cor}

In the subsequent sections, 
we shall give the proof of 
Theorem \ref{main},
which consists in the following steps:
\begin{enumerate}
\item
Well-definedness of the map $\Phi$.
\item
Well-definedness of the map $\Psi$.
\item
Bijectivity of $\Phi$ and $\Psi=\Phi^{-1}$.
\end{enumerate}

\renewcommand{\thesection}{\arabic{section}}
\section{Well-definedness of $\Phi$}
\setcounter{equation}{0}
\renewcommand{\theequation}{\thesection.\arabic{equation}}

For the well-definedness of $\Phi$,
it suffices to prove the following:
\begin{pro}\label{prop:phi}
Let $\lambda,\mu$ and $\nu$ be a Young diagram as before.
Suppose $f \in \PicA$. 
\begin{enumerate}
\item
The image $\Phi(f)$ is a Young 
tableau of shape $\mu$.
\item
Writing $R_{A'}(\Phi(f)) = \MECry$, the diagram 
$\lambda[i_1, i_2, \cdots ,i_k]$ is a Young diagram 
for any $k=1,\cd,N$ and $\lm[i_1,i_2,\cd,i_N]=\nu$.
\end{enumerate}
\end{pro}
\subsection{Proof of Proposition~\ref{prop:phi} (i)}

It is clear from the definition of $\Phi$ that $\Phi(f)$ 
is of shape $\mu$.
In order to prove (i), we may show
for any $i,j$:

{\bf (a)} $\Phi(f)_{i,j} <\Phi(f)_{i+1,j}$. \quad 
{\bf (b)} $\Phi(f)_{i,j} \leq \Phi(f)_{i,j+1}$.

\nd
{\bf (a)} By the definition of $\Phi$, we have
$\Phi(f)_{i,j} = f_1(i,j)$ and 
$\Phi(f)_{i+1,j} = f_1(i+1,j)$.
Since $f$ is a picture, we have
\begin{equation}
f(i,j) \leqslant_A f(i+1,j).\label{11a}
\end{equation}
Here assume  $f_1(i,j) \geq f_1(i+1,j)$.
If $f_2(i,j) \geq f_2(i+1,j)$, we have 
$f(i,j) \geqslant_P f(i+1,j)$.
Since $f$ is a picture, it means $(i,j)\geqslant_{A'} (i+1,j)$, 
which is a contradiction.
On the other hand, if $f_2(i,j) < f_2(i+1,j)$, 
then $f(i,j)\geqslant_A
f(i+1,j)$, which is also a contradiction.
Hence, we have $f_1(i,j) < f_1(i+1,j)$.

\nd
{\bf (b)}  As before, we have
$\Phi(f)_{i,j} = f_1(i,j)$ and $\Phi(f)_{i,j+1} = f_1(i,j+1)$.
Then we may show $f_1(i,j) < f_1(i,j+1)$, which will be shown 
 by the induction on $i$.
For the purpose, we need the following lemmas:
\begin{lem}\label{lem:1}
Let $f$ be in $\PicA$.
If $f_1(i,j) > f_1(i,j+1)$, then $f_2(i,j) > f_2(i,j+1)$.
\end{lem}
\textbf{Proof.}\quad
Suppose $f_2(i,j) \leq f_2(i,j+1)$.
Then by the assumption $f_1(i,j) > f_1(i,j+1)$, 
we obtain $f(i,j) \geqslant_A f(i,j+1)$.
However, since $(i,j) \leqslant_P (i,j+1)$ and $f$ is a picture, 
$f(i,j) \leqslant_A f(i,j+1)$, which derives a contradiction. Hence,
$f_2(i,j) > f_2(i,j+1)$.\qed

\begin{lem}\label{lem:2}
Suppose that $f_1(i,j) > f_1(i,j+1)$ for $f\in\PicA$.
Then there exists a unique $(k,l)$ in $\mu$ satisfying:
\begin{equation}
k<i,\,\,\; l \leq j,\q f_1(k,l)=f_1(i,j+1) \,\,\text{and}
\,\,f_2(k,l)=f_2(i,j).
\label{kl}
\end{equation}
\end{lem}
Note that by Lemma \ref{lem:1}, we have $f_2(i,j)>f_2(i,j+1)$:\\
\vskip-10mm
\setlength{\unitlength}{1pt}
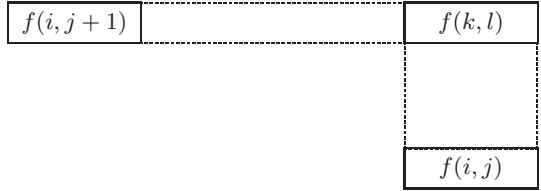
\begin{picture}(200,100)(0,0)
\put(50,55){\framebox(50,15){$f(i,j+1)$}}
\put(100,55){\dashbox(100,15){}}
\put(200,55){\framebox(50,15){$f(k,l)$}}
\put(200,15){\dashbox(50,40){}}
\put(200,0){\framebox(50,15){$f(i,j)$}}
\end{picture}

\textbf{Proof.}\quad
Since $\nu\setminus \lm$ is a skew diagram,
if $(a,b),\,\,(c,d)
\in\nu\setminus\lm$ satisfy $a<c,\,\,b<d$,
then $(a,d)\in\nu\setminus\lm$. Therefore, one gets
$(f_1(i,j+1),f_2(i,j))\in\nu\setminus\lm$.
It follows from the bijectivity of $f$ that there exists
a unique $(k,l)\in\mu$ such that
$f(k,l)=(f_1(i,j+1),f_2(i,j))$.
Now, it remains to show $k<i$ and $l \leq j$.
Since $(i,j) \leqslant_P (i,j+1)$, we have $f(i,j) \leqslant_A f(i,j+1)$ and then
\begin{equation}
f(k,l) \leqslant_A f(i,j) \leqslant_A f(i,j+1).
\label{klij}
\end{equation}
In the meanwhile, since
$f(i,j+1) \leqslant_P f(k,l) \leqslant_P f(i,j)$, we obtain
\begin{equation}
(i,j+1) \leqslant_{A'} (k,l) \leqslant_{A'} (i,j).
\label{ijkl}
\end{equation}
So, there can be $(k,l)$ in $X$ or $Y$
in the following figure, where
\[
X=\{(k,l)\in\mu \; : \; i<k \; \text{and} \;j<l\} \quad
Y=\{(k,l)\in\mu\;:\; i>k \;\text{and} \;l\leq j \}.
\]
\setlength{\unitlength}{1pt}
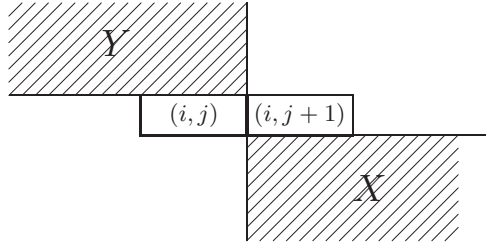
\begin{picture}(400,100)(0,0)
\put(140,50){\line(0,50){50}}
\put(100,65){\line(-50,0){50}}
\put(100,50){\framebox(40,15){$(i,j)$}}
\put(140,50){\framebox(40,15){$(i,j+1)$}}
\put(85,80){\Large{$Y$}}
\put(180,50){\line(50,0){50}}
\put(140,50){\line(0,-50){40}}
\put(180,25){\Large{$X$}}

\hatchpattern{5}{\path}
    \hatch(50,65)(140,100)
\hatchpattern{5}{\path}
    \hatch(140,10)(220,50)
\end{picture}
\\
In case  $(k,l) \in X$. 
Since $(i,j) \leqslant_P (k,l)$, we have
$f(i,j) \leqslant_A f(k,l)$. This does not match (\ref{klij}).
Hence $(k,l)\in Y$, i.e. $k<i$ and $l\leq j$.  \qed 

\vspace{0.3cm}

Let us show $f_1(i,j) \leq f_1(i,j+1)$ by the induction on $i$.

In case $i=1$, suppose $f_1(1,j) > f_1(1,j+1)$.
By Lemma~\ref{lem:1}, we have $f_2(1,j) > f_2(1,j+1)$.
By Lemma~\ref{lem:2},
there is $f(k,l)$ satisfying (\ref{kl}).
Since $i = 1$, the set $Y$ as above is, indeed, empty.
So there can not exist $(k,l)$. It contradicts  
Lemma~\ref{lem:2}. Thus, we obtain $f_1(1,j) \leq f_1(1,j+1)$.

In case $i=a>1$, assume $f_1(b,j)
\leq f_1(b,j+1)$ for any $b \leq a-1$ and any $j$.
If $f_1(a,j) > f_1(a,j+1)$,
by Lemma~\ref{lem:1} we have $f_2(a,j) > f_2(a,j+1)$.
In addition, by Lemma~\ref{lem:2},
there exists $f(k,l)$ satisfying (\ref{kl}) for $i=a$.
It follows from the hypothesis of the induction that
\[f_1(k,l)\leq f_1(k,l+1)\leq\cd \leq f_1(k,j+1),
\]
and we have
$f_1(k,j+1) <f_1(i,j+1)$ by (a).
Hence 
\[
f_1(k,l) < f_1(i,j+1),
\]
which contradicts (\ref{kl}). Hence,
we have $f_1(a,j) \leq f_1(a,j+1)$.
Now, we complete the proof of 
Proposition~\ref{prop:phi}(i).\qed
\subsection{Proof of Proposition~\ref{prop:phi} (ii)}

First, we prepare the following lemma:
\begin{lem}\label{lem:pic-no-ikisaki}
Let $f:\mu \to \nu \backslash \lm$ be a picture and
set
$ R_{A'}(\Phi(f)) = \MECry$.
Let $(p_j,q_j)\in\mu$ be the coordinate of
$\fsquare(0.4cm,i_j)$ in $\Phi(f)\in B(\mu)$ and
$(a_j, b_j) \in \nu$ the coordinate of the $j$-th
addition in $\lambda[i_1,\cdots,i_N]$. Then we have
$f(p_j,q_j) = (a_j, b_j)$ for any $j$.
\end{lem}
 
\textbf{Proof.}\quad
For any $m\in\{1,\cd, n+1\}$
list the coordinates in $\Phi(f)^{(m)}$ from the right as:
\[
(s_1,t_1), (s_2,t_2), \cdots,(s_c,t_c),
\]
where $c:=|\Phi(f)^{(m)}|$,
$s_1 \leq s_2 \leq \cdots \leq s_c$ and
$t_1 > t_2 > \cdots > t_c.$
So, we have
\begin{equation}
(s_1,t_1) \leqslant_{A'} (s_2,t_2) \leqslant_{A'} \cdots \leqslant_{A'} (s_c,t_c).
\label{stc}
\end{equation}
In the addition of $i_1,\cd,i_N$ to $\lm$, the box
$(s_k, t_k)$ goes to $(m, \lambda_m+k)\in \nu=\lm[i_1,\cd,i_N]$.

Here write the $m$-th row in $\nu \backslash \lambda$:
\[
(m,\lambda_m+1) \leqslant_P (m,\lambda_m+2) \leqslant_P \cdots \leqslant_P 
(m,\lambda_m+c)=(m,\nu_m).
\]
By the definition of $\Phi$, we have
\begin{equation}
 f^{-1}(\{(m,\lambda_m+1),\cd,(m,\lambda_m+c)\})=\Phi(f)^{(m)}
\label{phim}
\end{equation}
Since $f$ is a picture, we obtain
\begin{equation}
f^{-1}(m,\lambda_m+1) \leqslant_{A'} f^{-1}(m,\lambda_m+2) \leqslant_{A '}
\cdots \leqslant_{A'} f^{-1}(m,\lambda_m+c)=f^{-1}(m,\nu_m).
\label{fm}
\end{equation}
Thus, it follows from (\ref{stc}), (\ref{phim}) and 
(\ref{fm}) that 
$f^{-1}(m,\lambda_m+k) = (s_k,t_k)$ and then 
$f(s_k,t_k) = (m,\lambda_m+k)$ for any $k=1,\cd,c$ and 
$m=1\cd,n+1$.
\qed

\vspace{0.3cm}

\textbf{Proof of Proposition~\ref{prop:phi} (ii)}

Let $A'$ be an admissible order on $\mu$. 
We take the admissible reading $R_{A'}$ associated with $A'$
and write $R_{A'}(\Phi(f))=
\fsquare(0.3cm,i_1) \otimes \fsquare(0.3cm,i_2) 
\otimes \cdots \otimes\fsquare(0.3cm,i_N)
\in\mathbf{B}^{\otimes N}$.
Let us denote the coordinate of $\fsquare(0.3cm,i_k)$
by $(x_k,y_k)\in \bbN\times\bbN$.

We shall show that  $\lambda[i_1,i_2,\cdots,i_k]$ is 
a Young diagram for any $k$ by using 
the induction on $k$.
In case $k=1$, by the definition of the admissible reading,
$(x_1,y_1)$ is the minimum 
with respect to the order $\leqslant_{A'}$. Thus, we have
$(x_1,y_1) = (1,\mu_1)$.
Since $f$ is an admissible picture,
\begin{equation}
\hbox{$f(x_1,y_1)$ is minimal with respect to the 
order $\leqslant_P$ in $\nu  \backslash \lambda$.}
\label{minimal}
\end{equation}
Due to the definition of $\Phi$, 
we have $f_1(x_1,y_1)=i_1$.

Since $(1,\mu_1) = (x_1,y_1)$ is the right-most in $\mu$, 
the first entry $\Fsquare(0.4cm,i_1)$ goes to 
$(i_1,\lambda_{i_1}+1)$ by the addition. 
It follows from Lemma~\ref{lem:pic-no-ikisaki} that 
we have $f(x_1,y_1) = (i_1,\lambda_{i_1}+1)$.
Since $\lambda=(\lm_1,\lm_2,\cd)$ is a Young diagram, 
we obtain $\lambda_{i_1-1} \geq \lambda_{i_1}$.
Now, suppose $i_1>1$ and $\lambda_{i_1-1} = \lambda_{i_1}$.
Since we have $\lambda[i_1,i_2,\cdots,i_N]=\nu$, 
there exists some entry $\fsquare(0.3cm,i_k)$ 
$(k\geq 2)$ added to the coordinate 
$(i_1-1,\lambda_{i_1-1}+1)$, which 
means $f(x_k,y_k)=(i_1-1,\lambda_{i_1-1}+1)\leqslant_P
(i_1,\lambda_{i_1}+1)=f(x_1,y_1)$. Then it
contradicts the fact that  $f(x_1,y_1)$ is 
minimal with respect to the order $\leqslant_P$.
Thus, we have $i_1=1$ or $\lambda_{i_1-1} >\lambda_{i_1}$.
Then, $\lambda[i_1]$ is a Young diagram.

In case $k \geq 2$, assume that 
 $\lambda' = \lambda[i_1,i_2,\cdots,i_{k-1}]$ 
is a Young diagram.
The coordinate $(x_k,y_k)$ of $\fsquare(0.4cm,i_k)$ in $\mu$
is the minimum in $\mu\setminus 
\{(x_1,y_1),\cd,(x_{k-1},y_{k-1})\}$ 
with respect to the order $\leqslant_{A'}$. 
Since $f$ is an admissible picture, 
$f(x_k,y_k)$ is minimal in $\nu \backslash\lm'$ with respect to 
$\leqslant_P$. By the definition of $\Phi$, we have $f_1(x_k,y_k)=i_k$.

If $i_k=1$, trivially $\lm'[i_k]$ is a Young diagram.
In case $i_k >1$, by Lemma~\ref{lem:pic-no-ikisaki} we have
$f(x_k,y_k) = (i_k,\lambda'_{i_k}+1)$. 
Here, by the similar argument as above, we have 
$\lambda'_{i_k-1}> \lambda'_{i_k}$ and then 
$\lm'[i_k]$ is a Young diagram. Thus, 
$\lambda[i_1,i_2,\cdots,i_k]$ is a Young diagram
for any $k$.
It is trivial that $\lambda[i_1,i_2,\cdots,i_N]=\mu$ by the 
definition of $\Phi(f)$.
\qed

\section{Well-definedness of $\Psi$}

\begin{pro}\label{prop:psi}
For any $T$ in $\CryA$, 
we have $\Psi(T) \in \PicA$, that is, 
\begin{enumerate}
\item[\rm(1)]
$\Psi(T)$ is a map from $\mu$ to $\nu \backslash \lambda$.
\item[\rm(2)]
$\Psi(T)$ is a bijection.
\item[\rm(3)]
$\Psi(T)$ is an $(A,A')$-admissible picture.
\end{enumerate}
\end{pro}
The following lemma is needed to show Proposition~\ref{prop:psi}.
\begin{lem} \label{lem:tab-to-pic}
Let $T$ be in $\CryA$ and 
set $(a,b):=\Psi(T)(i,j)$ for $(i,j)$ in $\mu$.
Then the destination of $(i,j)$ by 
the addition of $R_{A'}(T)$ coincides with $(a,b)$.
\end{lem}

\textbf{Proof.}\quad
Set  $m := T_{i,j}$($=(i,j)$-entry in $T$).
Suppose that $(i,j)$ is the $p$-th element 
in $T^{(m)}$ from the right.
Then, by the definition of $\Psi$, we have
$\Psi(T) (i,j) = (m, \lambda_m +p)$.
On the other hand, in the course of the addition of $R_{A'}(T)$
to $\lm$, we see that $m=T_{i,j}$ is added $p$-th 
to the $m$-th row.
This means that $T_{i.j}$ goes to $(m,\lambda_m+p)$ 
by the addition.
This proves our claim. \qed
\subsection{Proof of Proposition~\ref{prop:psi}(1),(2)}
By the definition of $\Psi$, it is clear that 
$\Psi(T)$ is a map from $\mu$.
Let us denote $R_{A'}(T)=
\fsquare(4mm,i_1)\ot\cd\ot\fsquare(4mm,i_N)$.
Since $T \in \CryA$,
by the addition of $R_{A'}(T)$ to $\lambda$, we have
$\lm[i_1,\cd,i_N]=\nu$.
So the set of all the destinations by the addition 
coincides with $\nu \backslash \lambda$.
Then, by Lemma~\ref{lem:tab-to-pic}, it implies 
that $\Psi(T)(\mu) = \nu \backslash \lambda$. 
This also shows that $\Psi(T)$ is surjective and then
bijective, since $|\mu|=|\nu\setminus\lm|$.
\qed
%
%

\subsection{Proof of Proposition~\ref{prop:psi} (3)}
By the above results, we know that $f=\Psi(T)$ is a bijection.
Now, we shall prove
\begin{enumerate}
\item $f^{-1}$ is $PA'$-standard.
\item $f$ is $PA$-standard.
\end{enumerate}

Let us see (i).
Set $
f^{-1}(a,b) = (i,j)$, $f^{-1}(a,b+1) = (x,y)$, 
and $f^{-1}(a+1,b) = (s,t)$.
Thus, there exist $p$ , $q$ and $r$ such that
\begin{gather}
(a,b) = f(i,j) = (T_{i,j}, \lambda_{T_{i,j}} + p),
\label{ab}\\
(a,b+1) = f(x,y) = (T_{x,y}, \lambda_{T_{x,y}} + q),
\label{ab1}\\
(a+1,b) = f(s,t) = (T_{s,t}, \lambda_{T_{s,t}} + r).
\label{a1b}
\end{gather}
First, we shall show:
\[
f^{-1}(a,b) \leqslant_{A'} f^{-1}(a,b+1)
\quad \text{i.e.}\;(i,j) \leqslant_{A'} (x,y).
\]
Since $T_{i,j} = a = T_{x,y}$ by (\ref{ab}) and (\ref{ab1}),
we have $\lambda_{T_{i,j}} = \lambda_{T_{x,y}}$.
Furthermore, since
$b=\lambda_{T_{i,j}}+p, \; 
b+1=\lambda_{T_{x,y}}+q$ by (\ref{ab1}) and (\ref{a1b}), 
we obtain $q=p+1$.
The coordinates $(i,j)$ and $(x,y)$ are in $T^{(a)}$.
Therefore, $q=p+1$ means that $(i,j)$ 
is right to $(x,y)$ and 
then $(i,j) \leqslant_{A'} (x,y)$.

Next, let us see:
\[
 f^{-1}(a,b) \leqslant_{A'} f^{-1}(a+1,b)
\quad \text{i.e.} \quad(i,j) \leqslant_{A'} (s,t).
\]
Suppose $i \geq s$. From (\ref{ab}) and (\ref{a1b}), 
we have $T_{i,j} =a,T_{s,t}=a+1$.
So we obtain $T_{i.j} < T_{s,t}$, which means $j<t$.
This is as follows:

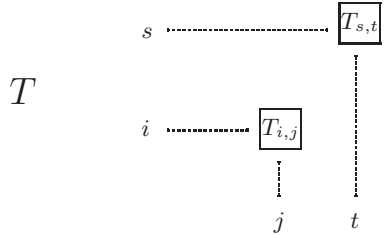
\setlength{\unitlength}{1pt}
\begin{picture}(400,100)(0,0)
\put(50,50){\Large{$T$}}
\put(100,75){$s$}
\put(110,78){\dashbox(60,0){}}
\put(175,73){\framebox(15,15){$T_{s,t}$}}
\put(100,38){$i$}
\put(110,40){\dashbox(30,0){}}
\put(145,33){\framebox(15,15){$T_{i,j}$}}
\put(150,3){$j$}
\put(152,15){\dashbox(0,13){}}
\put(179,3){$t$}
\put(181,15){\dashbox(0,53){}}
\end{picture}

Then we obtain $(i,j)\geqslant_{A'} (s,t)$.
In the process of the addition of $R_{A'}(T)
=\fsquare(4mm,i_1)\ot\cd\ot\fsquare(4mm,i_N)$ to $\lm$,
the entry at $(s,t)$ is added to $\lambda$ 
earlier than the one at $(i,j)$, which signifies 
that the coordinate $(a+1,b)$  is 
filled  earlier than the one at $(a,b)$
by the addition by Lemma \ref{lem:tab-to-pic}.
This means that there exists $k<N$ such that 
$\lm[i_1,\cd,i_k]$ is not a Young diagram,
which contradicts $T \in \CryA$.
Hence we have $i<s$.

Suppose $j<t$.  Owing to (\ref{ab}) and (\ref{a1b}), 
we obtain $T_{i,j} +1 = T_{s,t}$:

\setlength{\unitlength}{1pt}
\begin{picture}(400,100)(0,0)
\put(50,50){\Large{$T$}}
\put(100,78){$i$}
\put(110,80){\dashbox(30,0){}}
\put(175,33){\framebox(15,15){$T_{s,t}$}}
\put(100,35){$s$}
\put(110,38){\dashbox(60,0){}}
\put(145,73){\framebox(15,15){$T_{i,j}$}}
\put(147,3){$j$}
\put(149,15){\dashbox(0,53){}}
\put(179,3){$t$}
\put(181,15){\dashbox(0,13){}}
\end{picture}

Since $T$ is a Young tableau, we have
$a=T_{i,j}<T_{s,j}\leq T_{s,t}=a+1$ and then 
$T_{s,j}=a+1$. Thus, we know that there is no entry between
$(i,j)$ and $(s,j)$ and then $s=i+1$.
This is described as follows:

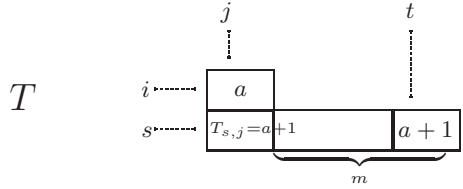
\setlength{\unitlength}{1pt}
\begin{picture}(400,100)(0,0)
\put(50,50){\Large{$T$}}

\put(100,55){$i$}
\put(105,58){\dashbox(15,0){}}
\put(125,50){\framebox(25,15){$a$}}
\put(130,85){$j$}
\put(133,70){\dashbox(0,10){}}

\put(100,40){$s$}
\put(105,43){\dashbox(15,0){}}
\put(125,35){\framebox(25,15){\q$\scriptstyle T_{s,j}=a+1$}}

\put(150,35){\framebox(45,15){}}
\put(195,35){\framebox(25,15){$a+1$}}
\put(200,85){$t$}
\put(202,55){\dashbox(0,25){}}

\put(150,35){$\underbrace{ \quad \quad \quad \quad \quad  \qquad}_{m}$}
\end{picture}

By the assumption $j< t$, we have $m:=t-j > 0$.
We have $f(s,j) = (a+1,\lambda_{a+1} +r+m)$:

\setlength{\unitlength}{1pt}
\begin{picture}(400,100)(0,0)
\put(50,50){\Large{$\nu\backslash\lambda$}}

\put(100,55){$a$}

\put(125,50){\framebox(28,15){$f(i,j)$}}
\put(135,20){$b$}

\put(100,40){$a+1$}
\put(125,35){\framebox(28,15){$f(s,t)$}}

\put(153,35){\framebox(40,15){}}
\put(193,35){\framebox(30,15){$f(s,j)$}}

\put(153,50){\dashbox(70,15){$\mathbf{X}$}}
\put(223,50){\dashbox(70,15){$\mathbf{Y}$}}
\put(223,35){\dashbox(70,15){$\mathbf{Z}$}}
\put(153,35){$\underbrace{\,\,\,
 \quad \quad \quad \quad \quad  \qquad}_{m}$}
\put(223,35){$\underbrace{\,\,\,
 \quad \quad \quad \quad \quad  \qquad}_{m}$}
\put(320,50){\rm Figure 1}
\end{picture}

Since $T \in \CryA$, the part \textbf{X} in the above 
figure must be filled earlier than $f(s,j)$ and later than 
$f(i,j)$ by the addition of $R_{A'}(T)$ to $\lm$.
Thus,  by Lemma \ref{lem:tab-to-pic},
$(i,j) \leqslant_{A'} f^{-1}(\mathbf{X}) \leqslant_{A'} (s,j)$, 
which implies that 
$f^{-1}(\mathbf{X})$ are in the shaded part 
in Figure 2 below:

\setlength{\unitlength}{1pt}
\begin{picture}(400,110)(0,0)
\put(100,50){\Large{$\mu$}}

\put(200,50){\line(-1,0){50}}
\put(200,65){\line(0,1){40}}

\hatchpattern{5}{\path}%
    \hatch(150,50)(200,105)

\put(200,50){\framebox(40,15){$(i,j)$}}
\put(200,35){\framebox(40,15){$(s,j)$}}

\put(240,50){\line(1,0){50}}
\put(240,50){\line(0,-1){40}}

\hatchpattern{5}{\path}%
    \hatch(240,10)(290,50)
\put(320,50){\rm Figure 2}
\end{picture}

Moreover, by the definition of $\Psi$,
an entry in $f^{-1}(\mathbf{X})$ is equal to 
$a$. Since the entries in $f^{-1}(\mathbf{X})$
are added later than the entry $T_{i,j}=a$ at $(i,j)$ and 
$T$ is a Young tableau, 
$f^{-1}(\mathbf{X})$ must be included in the shaded 
part in the following figure. 

\setlength{\unitlength}{1pt}
\begin{picture}(400,80)(0,0)
\put(100,50){\Large{$\mu$}}
\put(200,65){\line(-1,0){50}}
\put(200,40){\line(0,-1){20}}

\hatchpattern{5}{\path}%
    \hatch(150,20)(200,65)

\put(200,50){\framebox(30,15){$(i,j)$}}
\put(200,35){\framebox(30,15){$(s,j)$}}
\put(320,50){\rm Figure 3}
\end{picture}

\vskip-5mm
Therefore, by Figure 2 and 3, $f^{-1}(\mathbf{X})$ must be 
same as the shaded part in the following figure:
\vskip5mm

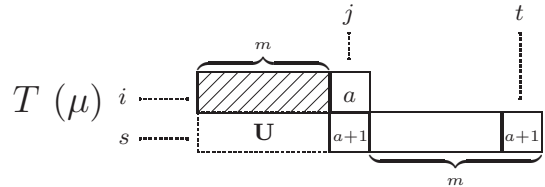
\setlength{\unitlength}{1pt}
\begin{picture}(400,80)(0,0)
\put(80,50){\Large{$T\,\,(\mu)$}}
\put(150,50){\framebox(50,15){}}

\hatchpattern{5}{\path}
    \hatch(150,50)(200,65)
\put(150,35){\dashbox(50,15){$\mathbf{U}$}}

\put(200,50){\framebox(15,15){}}
\put(200,35){\framebox(15,15){}}

\put(215,35){\framebox(50,15){}}
\put(265,35){\framebox(15,15){}}

\put(120,53){$i$}
\put(120,38){$s$}
\put(205,85){$j$}
\put(270,85){$t$}
\put(205,54){$a$}
\put(201,39){$\scriptstyle a+1$}
\put(266,39){$\scriptstyle a+1$}
\put(128,55){\dashbox(18,0){}}
\put(128,40){\dashbox(18,0){}}
\put(207,70){\dashbox(0,10){}}
\put(272,55){\dashbox(0,25){}}

\put(150,65){$\overbrace{\,\,  \quad \quad \quad \quad  \quad}^{m}$}
\put(215,35){$\underbrace{\,\,\,\,\quad\,\,
  \quad \quad \quad \quad  \quad}_{m}$}

\end{picture}

All the entries in this part are equal to $a$.
Then all the entries in the part $\mathbf U$
are equal to $a+1$.
Thus, we know that $f(\mathbf U)=\mathbf Z$
where $\mathbf Z$ is the part of $\nu\setminus\lm$ in Figure 1
as above. Arguing similarly, 
the part $\mathbf Y$ is sent to the part $\mathbf V$
in Figure 4 below and all the entries of $T$
in the part $\mathbf V$ are $a$. Then all the entries of $T$
in the part $\mathbf W$ are $a+1$.
The part $\mathbf W$ is sent to 
the right side of $\mathbf Z$ by $f$.

\setlength{\unitlength}{1pt}
\begin{picture}(400,80)(0,0)
\put(50,50){\Large{$T\,\,(\mu)$}}
\put(150,50){\framebox(50,15){}}
\put(165,56){$f^{-1}(\mathbf{X})$}
\put(150,35){\dashbox(50,15){$\mathbf{U}$}}
\put(100,35){\dashbox(50,15){$\mathbf{W}$}}
\put(100,50){\dashbox(50,15){$\mathbf{V}$}}

\put(200,50){\framebox(15,15){}}
\put(200,35){\framebox(15,15){}}

\put(215,35){\framebox(50,15){}}
\put(265,35){\framebox(15,15){}}

\put(205,54){$a$}
\put(201,39){$\scriptstyle a+1$}
\put(266,39){$\scriptstyle a+1$}

\put(150,65){$\overbrace{\,\,
  \quad \quad \quad \quad  \quad}^{m}$}
\put(215,35){$\underbrace{\,\,\,\,\quad\,\,
  \quad \quad \quad \quad  \quad}_{m}$}
\put(320,50){\rm Figure 4}

\end{picture}

Under the assumption 
$m=t-j>0$, we can repeat this process infinitely many times
and extend the $i$-th row and the $i+1$-row of $\mu$
unlimitedly,
which contradicts the 
finiteness of $\mu$. 
Thus, we have $ j \geq t$.
Finally, we have $i < s$ and $j \geq t$, and then 
$(i,j) \leqslant_{A'} (s,t)$. 
Then $f^{-1}$ is $PA'$-standard.\qed

Let us show (ii).
By the definition of $\Psi$, for $(i,j), (i,j+1),
(i+1,j)\in\mu$ there exist $p,q,$ and $r$ such 
that
\begin{gather}
f(i,j) = (T_{i,j},\lambda_{T_{i,j}} +p ),\label{ij}\\
f(i,j+1) = (T_{i,j+1},\lambda_{T_{i,j+1}} +q ),\label{ij1}\\
f(i+1,j) = (T_{i+1,j},\lambda_{T_{i+1,j}} +r ).\label{i1j}
\end{gather}
First, let us show 
\[
 f(i,j) \leqslant_A f(i,j+1).
\]
Since $T$ is a Young tableau, 
we have $T_{i,j} \leq T_{i,j+1}$ and then
 $f_1(i,j)=T_{i,j} \leq T_{i,j+1}=f_1(i,j+1)$, 
where $f=(f_1,f_2)$.

 In case $T_{i,j} = T_{i,j+1}$, 
we have $\lambda_{T_{i,j}} = \lambda_{T_{i,j+1}}$.
Moreover 
 $(i,j)$ is on the left-side of 
 $(i,j+1)$, 
which shows $p > q$.
Then, we have $\lambda_{T_{i,j}}+p > \lambda_{T_{i,j+1}}+q$, 
that is, $f_2(i,j) > f_2(i,j+1)$.
Hence we obtain $f(i,j) \leqslant_A f(i,j+1)$.

In case $T_{i,j} < T_{i,j+1}$, 
we have $\lambda_{T_{i,j}} \geq \lambda_{T_{i,j+1}}$.
Since $(i,j)\geqslant_{A'}(i,j+1)$, in the addition of 
$R_{A'}(T)$, the entry at 
$(i,j+1)$ is added earlier than the one at $(i,j)$.
Let $\lambda'$(resp. $\lm''$) be the resulting 
Young diagram obtained by the addition up to $(i,j+1)$
(resp. $(i,j)$).
It follows from  Lemma~\ref{lem:tab-to-pic} that 
the destination of $(i,j+1)$ (resp. $(i,j)$) 
by the addition coincides with 
$f(i,j+1) = (T_{i,j+1},\lambda_{T_{i,j}}+q)$
(resp. $f(i,j) = (T_{i,j},\lambda_{T_{i,j}}+p)$).

Since in  the all steps of the addition the resulting 
diagrams are always Young diagrams, we have
\[
\lambda'_{T_{i,j}} \geq \lambda'_{T_{i,j+1}} = \lambda_{T_{i,j+1}}+q.
\]
The entry $T_{i,j}$ is the $p$-th element 
in $T^{(T_{i,j})}$ from the right. 
Then by Lemma \ref{lem:tab-to-pic}, 
\[
\lambda'_{T_{i,j}}<\lm''_{T_{i,j}}
= \lambda_{T_{i,j}} +p.
\]
Hence 
\[
\lambda_{T_{i,j+1}}+q \leq 
\lambda'_{T_{i,j}} < \lambda_{T_{i,j}}+p,
\]
which implies $f_2(i,j+1) < f_2(i,j)$. 
We obtain $f(i,j) \leqslant_A f(i,j+1)$.

Next, let us show
\[
 f(i,j) \leqslant_A f(i+1,j).
\]
We have $f_1(i,j)=T_{i,j}<T_{i+1,j}=f_1(i+1,j)$. Thus,
 it is enough to show 
$\lambda_{T_{i,j}}+p \geq \lambda_{T_{i+1,j}}+r$ since 
$f_2(i,j)=\lambda_{T_{i,j}}+p$ and 
$f_2(i+1,j)=\lambda_{T_{i+1,j}}+r$, where $f=(f_1,f_2)$.

Assuming $f_2(i,j) < f_2(i+1,j)$, one has: 

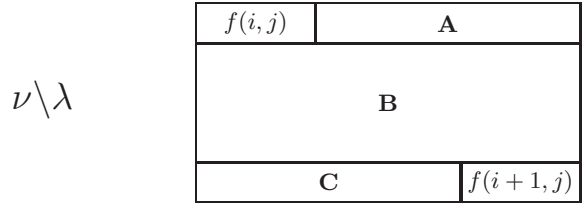
\setlength{\unitlength}{1pt}
\begin{picture}(400,100)(0,0)
\put(50,50){\Large{$\nu \backslash \lambda$}}
\put(120,75){\framebox(45,15){$f(i,j)$}}
\put(165,75){\framebox(100,15){$\mathbf{A}$}}
\put(120,30){\framebox(145,45){$\mathbf{B}$}}
\put(120,15){\framebox(100,15){$\mathbf{C}$}}
\put(220,15){\framebox(45,15){$f(i+1,j)$}}
\end{picture}
\begin{tabbing}
where \=${\mathbf A}=\{ (k,l) \; : \; 
k=f_1(i,j),\,\, f_2(i,j) < l \leq f_2(i+1,j) \}$,\\
\>${\mathbf B}=\{ (k,l) \; : \; 
f_1(i,j) <k<f_1(i+1,j),\,\; f_2(i,j) 
\leq l \leq f_2(i+1,j)\}$,\\
\>${\mathbf C}=\{ (k,l) \; : \; 
k=f_1(i+1,j),\,\, f_2(i,j) \leq l < f_2(i+1,j)\}$.
\end{tabbing}
Since $f^{-1}$ is $PA'$-standard and 
$f(i,j) \leqslant_P A,B,C \leqslant_P f(i+1,j)$,
we obtain $(i,j) \leqslant_{A'} f^{-1}(\mathbf A), 
f^{-1}(\mathbf B),$
$f^{-1}(\mathbf C) \leqslant_{A'} (i+1,j)$.
So the parts $f^{-1}(A), f^{-1}(B)$ and $ f^{-1}(C)$
must be in the shaded part of Figure 5 below:.

\setlength{\unitlength}{1pt}
\begin{picture}(400,100)(0,0)
\put(100,50){\Large{$\mu$}}

\put(200,50){\line(-1,0){50}}
\put(200,65){\line(0,1){30}}

\hatchpattern{5}{\path}
    \hatch(150,50)(200,95)

\put(200,50){\framebox(40,15){$(i,j)$}}
\put(200,35){\framebox(40,15){$(i+1,j)$}}

\put(240,50){\line(1,0){50}}
\put(240,50){\line(0,-1){40}}

\hatchpattern{5}{\path}
    \hatch(240,10)(290,50)
\put(350,50){\rm Figure 5}
\end{picture}

For $(k,l) \in {\mathbf B}$, 
by the definition of 
$\mathbf B$ as above, we have 
$T_{i,j}=f_1(i,j) < k < f_1(i+1,j)=T_{i+1,j}$.
If we set $(s,t):=f^{-1}(k,l)$, then $(s,t)$ is in the 
shaded part in Figure 3. Hence, 
in the Young tableau $T$ we have
$k=f_1(s,t)=T_{s,t}\leq T_{i,j}=f_1(i,j)$ 
or $k=T_{s,t}\geq T_{i+1,j}=f_1(i+1,j)$, 
which derives a contradiction.
Thus, $\mathbf B$ is empty and then the figure turns to:

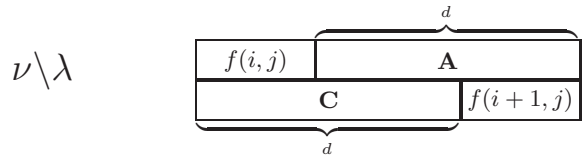
\setlength{\unitlength}{1pt}
\begin{picture}(400,100)(0,0)
\put(50,50){\Large{$\nu \backslash \lambda$}}
\put(120,50){\framebox(45,15){$f(i,j)$}}
\put(165,50){\framebox(100,15){$\mathbf{A}$}}

\put(120,35){\framebox(100,15){$\mathbf{C}$}}
\put(220,35){\framebox(45,15){$f(i+1,j)$}}

\put(165,65){$\overbrace{\,\,\,\, \quad \quad \quad \quad \quad  \qquad \qquad \quad}^{d}$}
\put(120,35){$\underbrace{\,\,\,\, \quad \quad \quad \quad \quad  \qquad\qquad \quad}_{d}$}
\end{picture}

\vskip-8mm
\noindent
where $d$ is the length of $\mathbf A$ and $\mathbf C$.
Now we have $T_{i,j}+1=f_1(i,j) +1 = f_1(i+1,j)=T_{i+1,j}$.
Since $T \in \CryA$, by the addition
the part $A$ is filled earlier than $f(i+1,j)$ is.
Since the part $\mathbf A$ is in the same 
row with $f(i,j)$, 
in the Young tableau $T$, the entries in the part 
$f^{-1}(\mathbf A)$ are equal to $T_{i,j}$. 
Then, there should be $f^{-1}(\mathbf A)\subset\mu$ 
in the shaded part of 
Figure 6 below:

\setlength{\unitlength}{1pt}
\begin{picture}(400,100)(0,0)
\put(100,50){\Large{$\mu$}}
\put(200,65){\line(-1,0){50}}
\put(200,40){\line(0,-1){20}}
\put(240,50){\line(1,0){50}}
\put(240,65){\line(0,1){25}}
\hatchpattern{5}{\path}
    \hatch(150,20)(200,65)
\hatchpattern{5}{\path}
    \hatch(240,50)(290,90)

\put(350,50){\rm Figure 6}
\put(200,50){\framebox(40,15){$(i,j)$}}
\put(200,35){\framebox(40,15){$(i+1,j)$}}

\end{picture}

By Figure 5 and 6, we know that 
$f^{-1}(A)$ is the following shaded part.

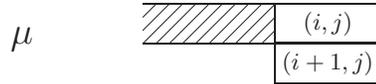
\setlength{\unitlength}{1pt}
\begin{picture}(400,70)(0,0)
\put(100,50){\Large{$\mu$}}
\put(200,65){\line(-1,0){50}}
\put(200,50){\line(-1,0){50}}

\hatchpattern{5}{\path}
    \hatch(150,50)(200,65)
\put(200,50){\framebox(40,15){$(i,j)$}}
\put(200,35){\framebox(40,15){$(i+1,j)$}}
\end{picture}

\vskip-8mm
Similar argument shows that
the entries in the part $\mathbf C$ are equal to 
$T_{i+1,j}=f_1(i+1,j)=T_{i,j}+1$ and the part 
$f^{-1}(C)$ is in the shaded part of the following figure.

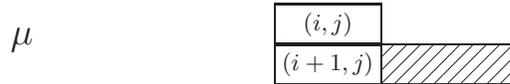
\setlength{\unitlength}{1pt}
\begin{picture}(400,80)(0,0)
\put(100,50){\Large{$\mu$}}
\put(240,50){\line(1,0){50}}
\put(240,35){\line(1,0){50}}
\hatchpattern{5}{\path}
    \hatch(240,35)(290,50)
\put(200,50){\framebox(40,15){$(i,j)$}}
\put(200,35){\framebox(40,15){$(i+1,j)$}}
\end{picture}

\vskip-8mm
Taking these into account, we get 

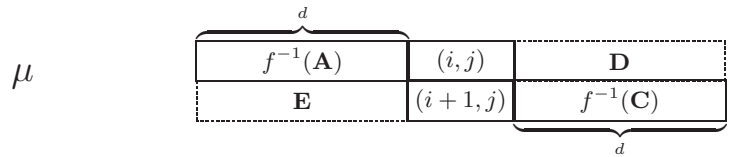
\setlength{\unitlength}{1pt}
\begin{picture}(400,80)(0,0)
\put(50,50){\Large{$\mu$}}
\put(120,50){\framebox(80,15){$f^{-1}(\mathbf{A})$}}

\put(200,50){\framebox(40,15){$(i,j)$}}
\put(200,35){\framebox(40,15){$(i+1,j)$}}

\put(240,50){\dashbox(80,15){$\mathbf{D}$}}
\put(240,35){\framebox(80,15){$f^{-1}(\mathbf{C})$}}
\put(120,35){\dashbox(80,15){$\mathbf{E}$}}

\put(120,65){$\overbrace{\,\,\,\, \quad \quad \quad \quad \quad  \qquad  \quad}^{d}$}
\put(240,35){$\underbrace{\,\,\,\, \quad \quad \quad \quad \quad  \qquad \quad}_{d}$}

\end{picture}

Since the difference between the entries 
$T_{i+1,j}$ and $T_{i,j}$ is just 1, 
the entries in the part $\mathbf D$ 
(resp. $\mathbf E$) above
are all equal to $T_{i,j}$ (resp. $T_{i+1,j}=T_{i,j}+1$).
Then we have the following figure:

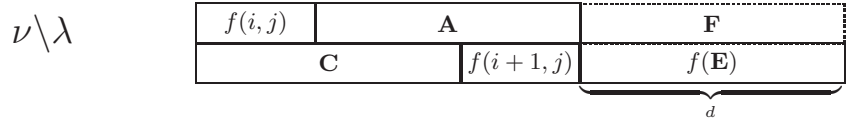
\setlength{\unitlength}{1pt}
\begin{picture}(400,100)(0,0)
\put(50,50){\Large{$\nu \backslash \lambda$}}
\put(120,50){\framebox(45,15){$f(i,j)$}}
\put(165,50){\framebox(100,15){$\mathbf{A}$}}
\put(265,50){\dashbox(100,15){$\mathbf{F}$}}
\put(120,35){\framebox(100,15){$\mathbf{C}$}}
\put(220,35){\framebox(45,15){$f(i+1,j)$}}
\put(265,35){\framebox(100,15){$f(\mathbf{E})$}}

\put(265,35){$\underbrace{\,\,\,\, \quad \quad 
\quad \quad \quad  \qquad  \qquad \quad}_{d}$}
\end{picture}

By the addition the part $\mathbf F$ should be filled 
earliear than $f(\mathbf E)$ is.
Then we obtain

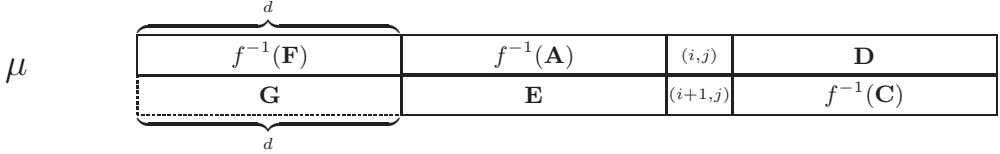
\setlength{\unitlength}{1pt}
\begin{picture}(500,100)(50,0)
\put(50,50){\Large{$\mu$}}
\put(100,50){\framebox(100,15){$f^{-1}(\mathbf{F})$}}
\put(200,50){\framebox(100,15){$f^{-1}(\mathbf{A})$}}
\put(300,50){\framebox(25,15){$\scriptstyle (i,j)$}}
\put(325,50){\framebox(100,15){$\mathbf{D}$}}

\put(100,35){\dashbox(100,15){$\mathbf{G}$}}
\put(200,35){\framebox(100,15){$\mathbf{E}$}}
\put(300,35){\framebox(25,15){$\scriptstyle (i+1,j)$}}
\put(325,35){\framebox(100,15){$f^{-1}(\mathbf{C})$}}

\put(100,65){$\overbrace{ \quad \quad \quad\,\,\,\,\quad \quad  \qquad \quad \quad \quad}^{d}$}
\put(100,35){$\underbrace{\quad \quad \quad\,\,\,\, \quad \quad  \qquad \quad \quad \quad}_{d}$}

\end{picture}
Arguing similarly, 
in the Young tableau $T$, all the entries in the part 
$f^{-1}(\mathbf F)$ above coincide with $T_{i,j}$ and then 
the ones in the part 
$\mathbf G$ coincide with $T_{i+1,j}=T_{i,j}+1$.
Then if $d>0$, repeating these arguments, we can 
extend the part in the left-side of $(i+1,j)$ 
(including the parts $\mathbf E$, $\mathbf G$)
unlimitedly. Indeed, it can not occur for the 
finite diagram $\mu$. Thus, we know that $d=0$.
This  means $\mathbf A=\mathbf C=\emptyset$, which 
contradicts the assumption $f_2(i,j)<f_2(i+1,j)$. 
Hence we have $f_2(i,j) \geq f_2(i+1,j)$, 
and then $f(i,j) \leqslant_A f(i+1,j)$.
It completes the proof.\qed

\section{Bijectivity of $\Phi$ and $\Psi$}

By the arguments above we get  the well-defined maps 
\[
\Phi\; :\; \PicA \longrightarrow \CryA
\qquad
\Psi \; :\; \CryA \longrightarrow \PicA.
\]
Now let us show:
\begin{enumerate}
\item $\Phi \circ \Psi = \text{id}_{\CryA}$
\item $\Psi \circ \Phi = \text{id}_{\PicA}$.
\end{enumerate}

(i) 
Recalling the definition of $\Psi$, 
for $T\in\CryA$ one gets that 
the admissible picture $\Psi(T)$ sends $(i,j)\in\mu$ to 
$(T_{i,j},  \lambda_{T_{i,j}}+p(i,j;T))$.
Due to the definition of $\Phi$ it is easy to see that 
$\Phi \circ \Psi(T)$ is an element in $\CryA$ 
whose $(i,j)$-entry is $T_{i,j}$, which implies 
$T=\Phi\circ\Psi(T)$.
Thus, we have $\Phi \circ \Psi=\text{id}_{\CryA}$ as desired.

(ii)
Taking an admissible  picture $f\in\PicA$, 
set $g:=\Psi \circ \Phi(f)\in\PicA$.
The image $\Phi(f)$ is an element in $\CryA$ whose $(i,j)$-entry is
$f_1(i,j)$. Then 
\begin{equation}
g(s,t) =(\Phi(f)_{s,t}, \lambda_{\Phi(f)_{s,t}}+p)
=(f_1(s,t), \lambda_{f_1(s,t)}+p),
\label{gst}
\end{equation}
where $(s,t)\in\mu$ is the $p$-th element from the right 
in $\Phi(f)^{(f_1(s,t))}:=\{(x,y)\in\mu|\Phi(f)_{x,y}=f_1(s,t)\}$.
It follows from  Lemma~\ref{lem:pic-no-ikisaki} that
\begin{equation}
f(s,t) = (\Phi(f)_{s,t} ,\lambda_{\Phi(f)_{s,t}}+p)
=(f_1(s,t), \lambda_{f_1(s,t)}+p).
\label{fst}
\end{equation}
which means $f(s,t) = g(s,t)$ and 
hence we have $\Psi \circ \Phi =$ id$_\PicA$.\qed

\end{document}